# A solution to the evolution-related Truscott-Brindley model for the generalized phytoplankton-zooplankton populations

**R. K. Upadhyay**

Department of Applied Mathematics
Indian School of Mines, Dhanbad, Jharkhand-826 004, India
E-mail: ranjit_ism@yahoo.com

**Abstract**

Phytoplankton are tiny floating plants (algae) living in oceans. In the process of photosynthesis, phytoplankton produces half of the world's oxygen. Moreover, by primary production, death and sinking, they transport carbon from the ocean's surface layer to marine sediments. There are many species of phytoplankton that can be distinguished according to morphology.

In this paper, we investigate the generalised Truscott-Brindley model of the dynamics of zoologically defined interacting populations which have spatial structure. Specifically, we consider conjointly marine phytoplankton and zooplankton populations, and model them as an excitable medium. The resolution using the *Boubaker polynomials expansion scheme* (*BPES*) along with stability analysis is carried out.

## 1. Introduction

The evolution-related Truscott-Brindley model [1] for the generalized phytoplankton-zooplankton populations is given by

$$\frac{dX}{dT} = rX\left(1 - \frac{X}{K}\right) - \omega Y \frac{X^2}{a^2 + X^2}, \tag{1a}$$

$$\frac{dY}{dT} = \beta \omega Y \frac{X^2}{a^2 + X^2} - \eta Y. \tag{1b}$$

where $X$ represents the population of phytoplankton and $Y$ the population of zooplankton. The growth rate of production of phytoplankton is represented by a logistic growth function, with a maximum growth rate $r$, and a carrying capacity $K$. Predation of the phytoplankton is represented by a Holling type III functional response [2], where $\omega$ is the maximum specific predation rate and $a$ governs how quickly that maximum is attained as prey density increases. Ecological observations [3] justify the use of a predation function which saturates for high prey densities. The parameter $\eta$ represents the specific rate of removal of zooplankton population. The parameter $\beta$ represents the ratio of biomass consumed to biomass of new herbivores produced. It covers a wide range of processes in nutrition and reproduction. A small proportion of the zooplankton are capable of reproduction, only some of the food ingested is assimilated and only a small amount of this is used for reproduction. It is not necessary to calculate these ratios directly as $\beta$ can be estimated implicitly from the values of steady population levels [3] and the value of maximum specific predation rate, $\omega$ [4]. The presence of $\beta$ or a similar effect is essential to the functioning of this model.

It is convenient to study the mathematical model given by equation (1) using its non-dimensional form. Using the substitutions

$$X = Kx, \ Y = Ky \ \text{and} \ T = t/\omega,$$

and denoting the non-dimensionless parameters as

$$\lambda = a/K, \ \alpha = r/\omega, \ \text{and} \ \gamma = \eta/\omega,$$

we obtain the following model system

$$\frac{dx}{dt} = \alpha\, x(1-x) - \frac{x^2 y}{\lambda^2 + x^2}, \tag{2a}$$

$$\frac{dy}{dt} = \frac{\beta\, x^2 y}{\lambda^2 + x^2} - \gamma\, y. \tag{2b}$$

The generalized Truscott- Brindley model, can be written as

$$\frac{dx}{dt} = \alpha\, x(t)(1-x(t)) - \left(\frac{(x(t))^{2+\mu}}{\lambda^2 + (x(t))^2}\right) y(t) \tag{3a}$$

$$\frac{dy}{dt} = \beta \left(\frac{(x(t))^{2+\mu}}{\lambda^2 + (x(t))^2}\right) y(t) - \gamma\, y(t) \tag{3b}$$

Stability is a main concern in population model modeling the predator-prey dynamics [25-27]. Dubey *et al.,* [17] proposed an analytical solution to the modified Lotka-Volterra model using Boubaker polynomial expansion scheme. Milgram commented on the claimed stability of solution to the accelerated-predator-satiety Lotka–Volterra predator–prey model, proposed by Dubey *et al.,*[17]. His critics are based on incompatibilities between the claimed asymptotic behavior and the presumed Malthusian growth of prey population in absence of predator. In the next section, I briefly explain the stability behavior of the proposed model system (3) where the growth of the prey population is governed by logistic growth function with a maximum growth rate $r$, and a carrying capacity $K$.

## 2. Stability Analysis

The model system (3a)-(3b) has three non-negative equilibrium points, $E_0 = (0,0)$, $E_1 = (1,0)$ and the intersection of the two isoclines is the equilibrium point $E^*(x^*, y^*)$, where $x^*$ and $y^*$ are related by $y^* = (\alpha\beta/\gamma) x^* (1-x^*)$. Since, $x^* \geq 0$ and $y^* \geq 0$, we have the condition $1 - x^* \geq 0$.

By computing the variational matrices corresponding to each equilibrium point, we obtain the following results:

(a) $E_0$ is a saddle point with unstable manifold locally in the x-direction and with stable manifold locally in the y-direction.

(b) $E_1$ is locally asymptotically stable if $\dfrac{1^{1+\mu}}{1+\lambda^2} < \dfrac{\gamma}{\beta}$ and is a saddle point if $\dfrac{1^{1+\mu}}{1+\lambda^2} > \dfrac{\gamma}{\beta}$.

(c) The interior equilibrium point $E^*$ is locally asymptotically stable in the $xy$-plane if and only if

$$\frac{x^*-1}{\lambda^2+x^{*2}} < \frac{x^*}{(\mu+1)\lambda^2+(\mu-1)x^{*2}}.$$

This result follows from the Routh–Hurwitz criteria.

The global stability behaviour of the interior equilibrium point $E^*(x^*, y^*)$ is given in the following theorem.

**Theorem** The equilibrium point $E^*(x^*, y^*)$ is globally asymptotically stable if the following condition holds:

$$\beta\{\lambda^2 P + x^2(x^*)^2 Q\} < x^{1+\mu}\{(x^*)^2 + \lambda^2)\} \tag{4}$$

where

$$P = x^{\mu+1} + x^\mu x^* + x^{\mu-1}(x^*)^2 + \ldots + x^2(x^*)^{\mu-1} + x(x^*)^\mu + (x^*)^{\mu+1}$$

and $\quad Q = x^{\mu-1} + x^{\mu-2}x^* + x^{\mu-3}(x^*)^2 + \ldots + x^2(x^*)^{\mu-3} + x(x^*)^{\mu-2} + (x^*)^{\mu-1}.$

**Proof** The proof is based on a Lyapunov direct method. Consider the positive definite function $V(x, y, z) = V_1(x, y, z) + V_2(x, y, z)$ where

$$V_1 = C_1\left[x - x^* - x^* \ln\left(\frac{x}{x^*}\right)\right]; \quad V_2 = C_2\left[y - y^* - y^* \ln\left(\frac{y}{y^*}\right)\right].$$

Now, V is a continuous function on $R_+^3 = \{(x, y, z) \in R^3 : x > 0, y > 0, z > 0\}$, and $C_1, C_2$ are positive constants to be determined. In order to investigate the global dynamics of the non-negative equilibrium point $E^*$ of the model system (3), the derivative of V with respect to time is computed as

$$\frac{dV}{dt} = \frac{dV_1}{dt} + \frac{dV_2}{dt}. \tag{5}$$

We have

$$\frac{dV_1}{dt} = -C_1\left[\alpha + \frac{\{x^2(x^*)^2 M + \lambda^2 N\}y^*}{(\lambda^2+x^2)\{\lambda^2+(x^*)^2\}}\right](x-x^*)^2 - C_1\left[\frac{x^{1+\mu}(x-x^*)(y-y^*)}{(\lambda^2+x^2)}\right] \tag{6a}$$

$$\frac{dV_2}{dt} = C_2\beta\left[\frac{\{\lambda^2 P + x^2(x^*)^2 Q\}}{(\lambda^2+x^2)\{\lambda^2+(x^*)^2\}}\right](x-x^*)(y-y^*), \tag{6b}$$

where

$$M = x^{\mu-2} + x^{\mu-3}x^* + x^{\mu-4}(x^*)^2 + \ldots + x^2(x^*)^{\mu-4} + x(x^*)^{\mu-3} + (x^*)^{\mu-2},$$

$$N = x^\mu + x^{\mu-1}x^* + x^{\mu-2}(x^*)^2 + \ldots + x^2(x^*)^{\mu-2} + x(x^*)^{\mu-1} + (x^*)^\mu,$$

$$P = x^{\mu+1} + x^\mu x^* + x^{\mu-1}(x^*)^2 + \ldots + x^2(x^*)^{\mu-1} + x(x^*)^\mu + (x^*)^{\mu+1},$$

$$Q = x^{\mu-1} + x^{\mu-2}x^* + x^{\mu-3}(x^*)^2 + \ldots + x^2(x^*)^{\mu-3} + x(x^*)^{\mu-2} + (x^*)^{\mu-1}.$$

Substituting (6) in (5) we obtain:

$$\frac{dV}{dt} = -C_1\left[\alpha + \frac{\{x^2(x^*)^2 M + \lambda^2 N\}y^*}{(\lambda^2 + x^2)\{\lambda^2 + (x^*)^2\}}\right](x - x^*)^2 
- \left[C_1 \frac{x^{1+\mu}}{(\lambda^2 + x^2)} - C_2\beta\left(\frac{\{\lambda^2 P + x^2(x^*)^2 Q\}}{(\lambda^2 + x^2)\{\lambda^2 + (x^*)^2\}}\right)\right](x - x^*)(y - y^*) \quad (7)$$

Without loss of generality, we choose $C_1 = C_2 = 1/2$.

Clearly, $dV/dt$ is negative under the condition (4). We conclude that $V$ is a Lyapunov function with respect to $E^*$. Hence, $E^*$ is globally asymptotically stable.

## 3. BPES resolution protocol

The resolution of system (3) is based on combining the two equations, as a prelude to the application of the Boubaker polynomials expansion scheme (BPES). The system can be written as

$$\begin{cases} \beta \dfrac{dx(t)}{dt} + \dfrac{dy(t)}{dt} = \beta\alpha\, x(t)(1 - x(t)) - \gamma\, y(t) \\ \dfrac{dy(t)}{dt} = \beta\left(\dfrac{(x(t))^{2+\mu}}{\lambda^2 + (x(t))^2}\right)y(t) - \gamma\, y(t) \end{cases} \quad (8)$$

The Boubaker polynomials expansion scheme (BPES) ([5-16]) is an analytical resolution protocol which was published by Oyodum et al. [5]. Recently, Slama et al. ([8-10]) developed a numerical model for the spatial time-dependant evolution of $A3$ melting point in $C40$ steel material during a particular sequence of resistance spot welding. In ecology, a solution of accelerated-predator-satiety Lotka-Volterra predator-prey model is obtained by Dubey et. al. [17] using the Boubaker polynomials expansion scheme.

In the present study, the BPES protocol investigations are performed using the following series expansion

$$x(t) = \frac{1}{2N_0} \sum_{q=1}^{N_0} \xi_q B_{4q}\left(v_q \frac{t}{t_m}\right) = \frac{1}{2N_0} \sum_{q=1}^{N_0} \xi_q B_{4q}(t^*), \quad x(0) = x_0,$$

$$y(t) = \frac{1}{2N_0} \sum_{q=1}^{N_0} \xi'_q B_{4q}(t^*), \quad t^* = v_q \frac{t}{t_m}, \quad y(0) = y_0. \quad (9)$$

where $v_q$ is the $4q$-Boubaker polynomial minimal root, $N_0$ is a pre-fixed integer, $t_m = 2\pi/\sqrt{\alpha\gamma}$ is the system characteristic time constant, and $\xi_q$ and $\xi'_q$, $q = 1, 2, ..., N_0$ are the unknown coefficients to be determined.

The initial conditions imply $\sum_{q=1}^{N_0} \xi_q = -N_0 x_0$, $\sum_{q=1}^{N_0} \xi'_q = -N_0 y_0$, since $B_{4q}(0) = -2$.

Substituting the expressions (9) in the first equation of (8), we obtain

$$\sum_{q=1}^{N_0} [\beta \xi_q + \xi'_q] \left(\frac{v_q}{t_m}\right) \frac{d}{dt^*} B_{4q}(t^*)$$

$$= \beta\alpha \left[ \left\{\sum_{q=1}^{N_0} \xi_q B_{4q}(t^*)\right\} \left\{1 - \frac{1}{2N_0} \sum_{q=1}^{N_0} \xi_q B_{4q}(t^*)\right\} \right] - \gamma \sum_{q=1}^{N_0} \xi'_q B_{4q}(t^*). \quad (10)$$

Denote $\quad W_q = \beta\left(\frac{v_q}{t_m}\right) \frac{d}{dt^*} B_{4q}(t^*) - \beta\alpha B_{4q}(t^*), \quad M_q = \sqrt{\frac{\beta\alpha}{2N_0}} B_{4q}(t^*),$

$$Z_q = -\left(\frac{v_q}{t_m}\right) \frac{d}{dt^*} B_{4q}(t^*) - \gamma B_{4q}(t^*). \quad t^* = v_q \frac{t}{t_m}.$$

Then, (10) can be written as

$$\sum_{q=1}^{N_0} \xi_q \left\{W_q + M_q \left[\sum_{q=1}^{N_0} \xi_q M_q\right]\right\} = \sum_{q=1}^{N_0} \xi'_q Z_q.$$

The non-zero solution is obtained by minimizing the $N_0$-dependent function $\Lambda_{N_0}$

$$\Lambda_{N_0} = \int_0^{t_m} \left[ \sum_{q=1}^{N_0} \xi_q \left\{W_q + M_q \left[\sum_{q=1}^{N_0} \xi_q M_q\right]\right\} - \sum_{q=1}^{N_0} \xi'_q Z_q \right] dt, \quad (11)$$

under the two initial conditions $\quad \sum_{q=1}^{N_0} \xi_q = -N_0 x_0, \quad \sum_{q=1}^{N_0} \xi'_q = -N_0 y_0.$

Since the non-zero solution is not unique, the same algorithm is applied to the second equation. It is equivalent to selecting a solution among the yielded solutions. For a given value of $\lambda$, the set of solutions that minimizes the $N_0$-dependent function $\Lambda'_{N_0}$:

$$\frac{dy(t)}{dt} = \beta \left( \frac{(x(t))^{2+\mu}}{\lambda^2 + (x(t))^2} \right) y(t) - \gamma\, y(t)$$

Denote $\quad P_q = \left(\dfrac{v_q}{t_m}\right) \dfrac{d}{dt^*} B_{4q}(t^*). \quad t^* = v_q \dfrac{t}{t_m}.$

$$\Lambda'_{N_0} = \int_0^{t_m} \left\{ \sum_{q=1}^{N_0} \xi'_q P_q - \left[ \frac{\beta \left\{ \frac{1}{2N_0} \sum_{q=1}^{N_0} \xi_q B_{4q}(t^*) \right\}^{2+\mu}}{\lambda^2 + \left\{ \frac{1}{2N_0} \sum_{q=1}^{N_0} \xi_q B_{4q}(t^*) \right\}^2} - \gamma \right] \sum_{q=1}^{N_0} \xi'_q B_{4q}(t^*) \right\} dt \quad (12)$$

The optimal values of $\xi_q$ and $\xi'_q$ may be generated by minimizing the sum of $\Lambda_{N_0}$ and $\Lambda'_{N_0}$ for some given value of $\lambda$ under the condition (12).

### 4. Numerical simulation

For the numerical study of the model system (3), we have considered two sets of parameter values. For the first set of parameters values, we obtain stable focus and limit cycles as presented in Figure 1 for different values of an adjustable non-negative parameter $\mu = 0, 0.5$ and 1.

For the first study, the initial conditions are taken as $x_0 = 0.9,\ y_0 = 0.5,$ and $t_m = 6.66$.

The following dynamical outcomes are obtained.

**Case 1:** $\mu = 0$:

Stable focus:

(i) $0.002 \le \alpha \le 5,\ \lambda = 0.057,\ \beta = 1.3,\ \gamma = 0.5$.
(ii) $\alpha = 1.9,\ \lambda = 0.003$ and $0.005 \le \lambda \le 1.26,\ \beta = 1.3,\ \gamma = 0.5$.
(iii) $\alpha = 1.9,\ \lambda = 0.057,\ \gamma = 0.5,\ 0.94 \le \beta \le 5.0$.
(iv) $\alpha = 1.9,\ \lambda = 0.057,\ \beta = 1.3,\ 0.001 \le \gamma \le 0.69,$ and $3.2 \le \gamma \le 5.0$.

Stable limit cycle:

(i) $\alpha = 1.9,\ \lambda = 0.057,\ \gamma = 0.5,\ 0.51 \le \beta \le 0.93$.
(ii) $\alpha = 1.9,\ \lambda = 0.057,\ \beta = 1.3,\ 0.7 \le \gamma \le 1.28$.

Extinction:

    (i) $\alpha = 0.001,\ \lambda = 0.057, \beta = 1.3, \gamma = 0.5$.
    (ii) $\alpha = 1.9, \lambda = 0.003$ and $1.27 \leq \lambda \leq 5.0,\ \beta = 1.3, \gamma = 0.5$.
    (iii) $\alpha = 1.9, \lambda = 0.057,\ \beta = 1.3, 1.3 \leq \gamma \leq 3.1$.

Integration error: (i) $\alpha = 1.9, 0.001 \leq \lambda \leq 0.002$, and $\lambda = 0.004,\ \beta = 1.3, \gamma = 0.5$.

**Case 2:** $\mu = 0.5$:

Stable focus:

    (i) $\alpha = 1.9, 0.076 \leq \lambda \leq 1.26,\ \beta = 1.3, \gamma = 0.5$'
    (ii) $\alpha = 1.9, \lambda = 0.057,\ \gamma = 0.5, 0.5 \leq \beta \leq 0.94, 1.76 \leq \beta \leq 2.29$.
    (iii) $\alpha = 1.9, \lambda = 0.057,\ \beta = 1.3, 0.005 \leq \gamma \leq 0.03, 0.21 \leq \gamma \leq 0.32$, and $0.69 \leq \gamma \leq 1.29$.

Stable limit cycle:

    (i) $\beta = 1.3, \gamma = 0.5,\ \lambda = 0.057,\ 0.1 \leq \alpha \leq 0.21, 0.31 \leq \alpha \leq 0.33, \alpha = 0.35$,
        $0.37 \leq \alpha \leq 0.4,\ 0.42 \leq \alpha \leq 0.44,\ 0.57 \leq \alpha \leq 0.85, 0.91 \leq \alpha \leq 5.0$.
    (ii) $\alpha = 1.9, 0.051 \leq \lambda \leq 0.075, \gamma = 0.5, \beta = 1.3$.
    (iii) $\alpha = 1.9, \lambda = 0.057, 0.95 \leq \beta \leq 1.75, \gamma = 0.5$.
    (iv) $\alpha = 1.9, \lambda = 0.057,\ \beta = 1.3, 0.41 \leq \gamma \leq 0.68$.

Extinction:

    (i) $\alpha = 1.9, 1.27 \leq \lambda \leq 5.0,\ \beta = 1.3, \gamma = 0.5$.
    (ii) $\alpha = 1.9, \lambda = 0.057,\ 0.001 \leq \beta \leq 0.49, \gamma = 0.5$.
    (iii) $\alpha = 1.9, \lambda = 0.057,\ \beta = 1.3, 1.3 \leq \gamma \leq 5.0$.

Argument domain error:

    (i) $\beta = 1.3, \gamma = 0.5,\ \lambda = 0.057,\ 0.001 \leq \alpha \leq 0.09, 0.22 \leq \alpha \leq 0.3, \alpha = 0.34, \alpha = 0.36$,
        $\alpha = 0.41, 0.45 \leq \alpha \leq 0.56,\ 0.86 \leq \alpha \leq 0.9$.
    (ii) $\alpha = 1.9, 0.001 \leq \lambda \leq 0.5,\ \beta = 1.3, \gamma = 0.5$.
    (iii) $\alpha = 1.9, \lambda = 0.057,\ 2.3 \leq \beta \leq 5.0, \gamma = 0.5$.
    (iv) $\alpha = 1.9, \lambda = 0.057,\ \beta = 1.3, 0.001 \leq \gamma \leq 0.004, 0.04 \leq \gamma \leq 0.2, 0.33 \leq \gamma \leq 0.4$.

**Case 3:** $\mu = 1.0$:

Stable focus:

    (i) $0.001 \leq \alpha \leq 5.0,\ \lambda = 0.057,\ \beta = 1.3, \gamma = 0.5$.
    (ii) $\alpha = 1.9, 0.001 \leq \lambda \leq 1.27,\ \beta = 1.3, \gamma = 0.5$.
    (iii) $\alpha = 1.9, \lambda = 0.057,\ \gamma = 0.5, 0.5 \leq \beta \leq 5.0$.

(iv) $\alpha = 1.9, \lambda = 0.057, \beta = 1.3, 0.001 \leq \gamma \leq 1.29$, and $3.2 \leq \gamma \leq 5.0$.

Extinction:

(i) $\alpha = 1.9, \lambda = 0.003$ and $1.28 \leq \lambda \leq 5.0$, $\beta = 1.3, \gamma = 0.5$.
(ii) $\alpha = 1.9, \lambda = 0.057, \beta = 1.3, 1.3 \leq \gamma \leq 3.19$.

Simulation is also done taking the initial values as $x_0 = 0.5$, $y_0 = 0.5$, and $t_m = 8.88$.

The following dynamical outcomes are obtained.

**Case 1:** $\mu = 0$:

Stable focus:

(i) $0.002 \leq \alpha \leq 5$, $\lambda = 0.057, \beta = 1.2, \gamma = 0.5$.
(ii) $\alpha = 1.0$, $0.006 \leq \lambda \leq 1.18$, $\beta = 1.2, \gamma = 0.5$.
(iii) $\alpha = 1.0, \lambda = 0.057, \gamma = 0.5, 0.94 \leq \beta \leq 5.0$.
(iv) $\alpha = 1.0, \lambda = 0.057, \beta = 1.2, 0.004 \leq \gamma \leq 0.63, \gamma = 1.19$, and $2.2 \leq \gamma \leq 5.0$.

Stable limit cycle:

(i) $\alpha = 1.0, \lambda = 0.057, \gamma = 0.5, 0.51 \leq \beta \leq 0.93$.
(ii) $\alpha = 1.0, \lambda = 0.057, \beta = 1.2, 0.64 \leq \gamma \leq 1.18$.

Extinction:

(i) $\alpha = 0.001$, $\lambda = 0.057, \beta = 1.2, \gamma = 0.5$.
(ii) $\alpha = 1.0, \lambda = 0.002$ and $1.19 \leq \lambda \leq 5.0$, $\beta = 1.2, \gamma = 0.5$.
(iii) $\alpha = 1.0, \lambda = 0.057, \gamma = 0.5, 0.001 \leq \beta \leq 0.5$.
(iv) $\alpha = 1.0, \lambda = 0.057, \beta = 1.2, 0.001 \leq \gamma \leq 0.003$, and $1.2 \leq \gamma \leq 2.19$.

Integration error:

(i) $\alpha = 1.0, \lambda = 0.001$, and $0.003 \leq \lambda \leq 0.005$, $\beta = 1.2, \gamma = 0.5$.

**Case 2:** $\mu = 0.5$:

Stable focus:

(i) $\alpha = 1.0, 0.078 \leq \lambda \leq 1.18$, $\beta = 1.2, \gamma = 0.5$.
(ii) $\alpha = 1.0, \lambda = 0.057, \gamma = 0.5, 0.51 \leq \beta \leq 0.94, 1.76 \leq \beta \leq 5.0$.
(iii) $\alpha = 1.0, \lambda = 0.057, \beta = 1.2, 0.001 \leq \gamma \leq 0.08, 0.12 \leq \gamma \leq 0.34$, $0.64 \leq \gamma \leq 1.19$, and $2.2 \leq \gamma \leq 5.0$.

Stable limit cycle:

  (i) $\beta = 1.2, \gamma = 0.5$, $\lambda = 0.057$, $0.001 \leq \alpha \leq 5.0$.
  (ii) $\alpha = 1.0, \gamma = 0.5, \beta = 1.2, 0.0036 \leq \lambda \leq 0.054, 0.056 \leq \lambda \leq 0.059, 0.061 \leq \lambda \leq 0.077$.
  (iii) $\alpha = 1.0, \lambda = 0.057, 0.95 \leq \beta \leq 1.75, \gamma = 0.5$.
  (iv) $\alpha = 1.0, \lambda = 0.057, \beta = 1.2, 0.35 \leq \gamma \leq 0.63$.

Extinction:

  (i) $\alpha = 1.0, 1.19 \leq \lambda \leq 5.0, \beta = 1.2, \gamma = 0.5$.
  (ii) $\alpha = 1.0, \lambda = 0.057$, $0.001 \leq \beta \leq 0.5, \gamma = 0.5$.
  (iii) $\alpha = 1.0, \lambda = 0.057, \beta = 1.2, 1.2 \leq \gamma \leq 2.19$.

Argument domain error:

  (i) $\alpha = 1.0, \beta = 1.2, \gamma = 0.5$, $0.001 \leq \lambda \leq 0.035, \lambda = 0.55, \lambda = 0.06$,
  (ii) $\alpha = 1.0, \lambda = 0.057, \beta = 1.2, 0.09 \leq \gamma \leq 0.11$.

**Case 3:** $\mu = 1.0$:

Stable focus:

  (i) $0.001 \leq \alpha \leq 5.0$, $\lambda = 0.057, \beta = 1.2, \gamma = 0.5$.
  (ii) $\alpha = 1.0, 0.001 \leq \lambda \leq 1.19$, $\beta = 1.2, \gamma = 0.5$.
  (iii) $\alpha = 1.0, \lambda = 0.057, \gamma = 0.5, 0.51 \leq \beta \leq 5.0$.
  (iv) $\alpha = 1.0, \lambda = 0.057, \beta = 1.2, 0.001 \leq \gamma \leq 1.19$, and $2.2 \leq \gamma \leq 5.0$.

Extinction:

  (i) $\alpha = 1.0, 0.002 \leq \lambda \leq 5.0, \beta = 1.2, \gamma = 0.5$.
  (ii) $\alpha = 1.0, \lambda = 0.057, \gamma = 0.5, 0.001 \leq \beta \leq 0.5$.
  (iii) $\alpha = 1.0, \lambda = 0.057, \beta = 1.2, 1.2 \leq \gamma \leq 2.19$.

For the parameters set, $\alpha = 1.9, \lambda = 0.057, \beta = 1.3, \gamma = 0.5, t_m = 6.66, x_0 = 0.9, y_0 = 0.5$, we have plotted the phase portrait and its time-series which are given below (Figs. 1a, 1b, 1c):

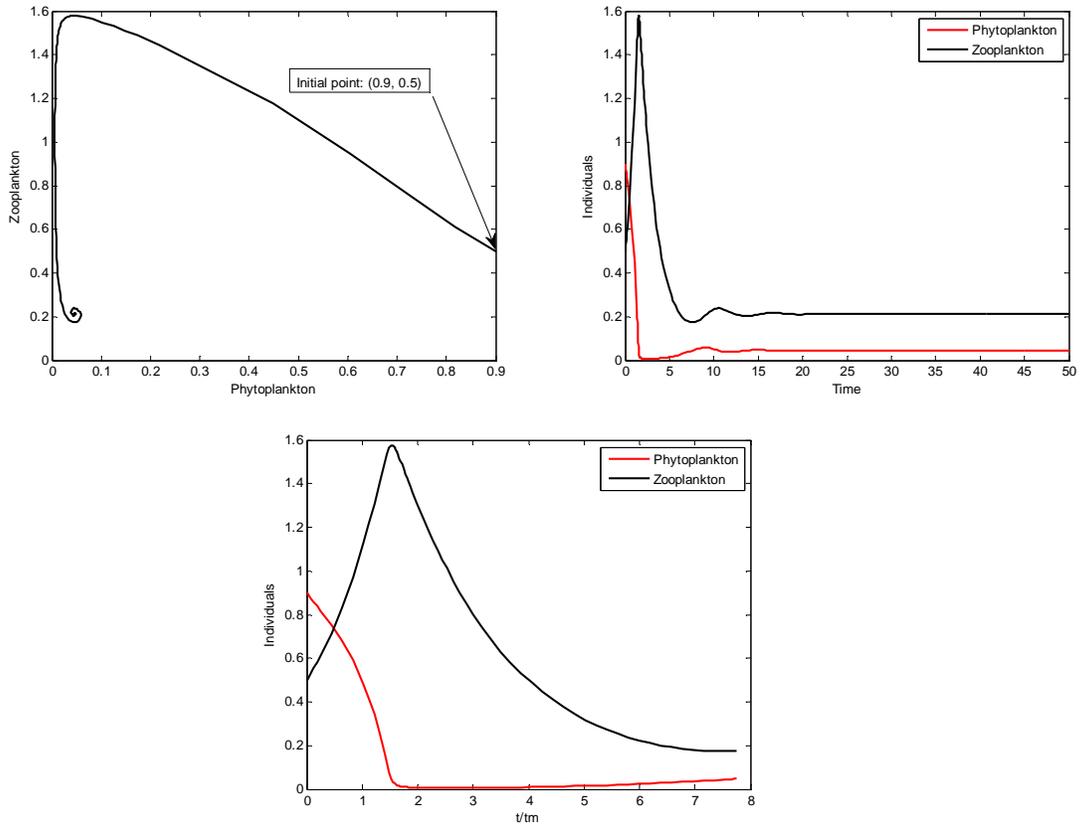

Fig. 1a. Phytoplankton-zooplankton evolution and orbit associated with $\mu = 0$.

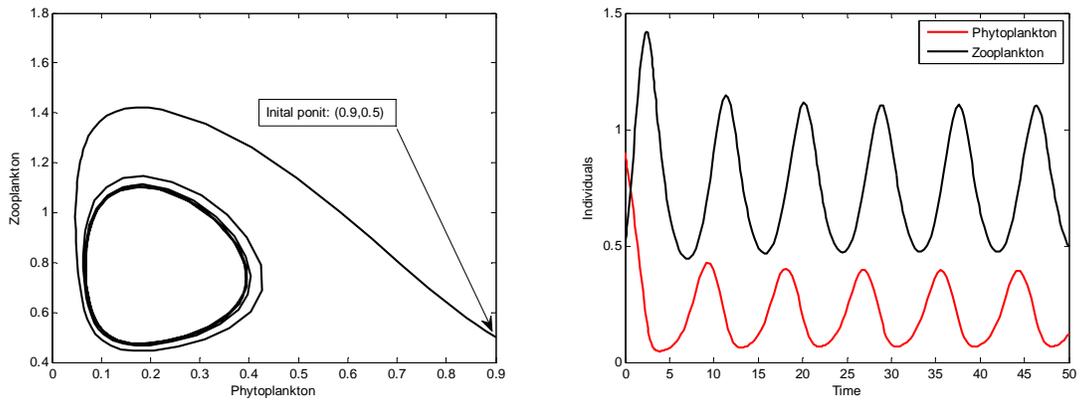

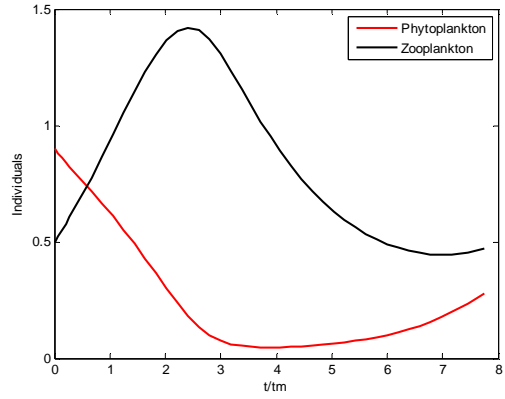

Fig. 1b. Phytoplankton-zooplankton evolution and orbit associated with $\mu = 0.5$.

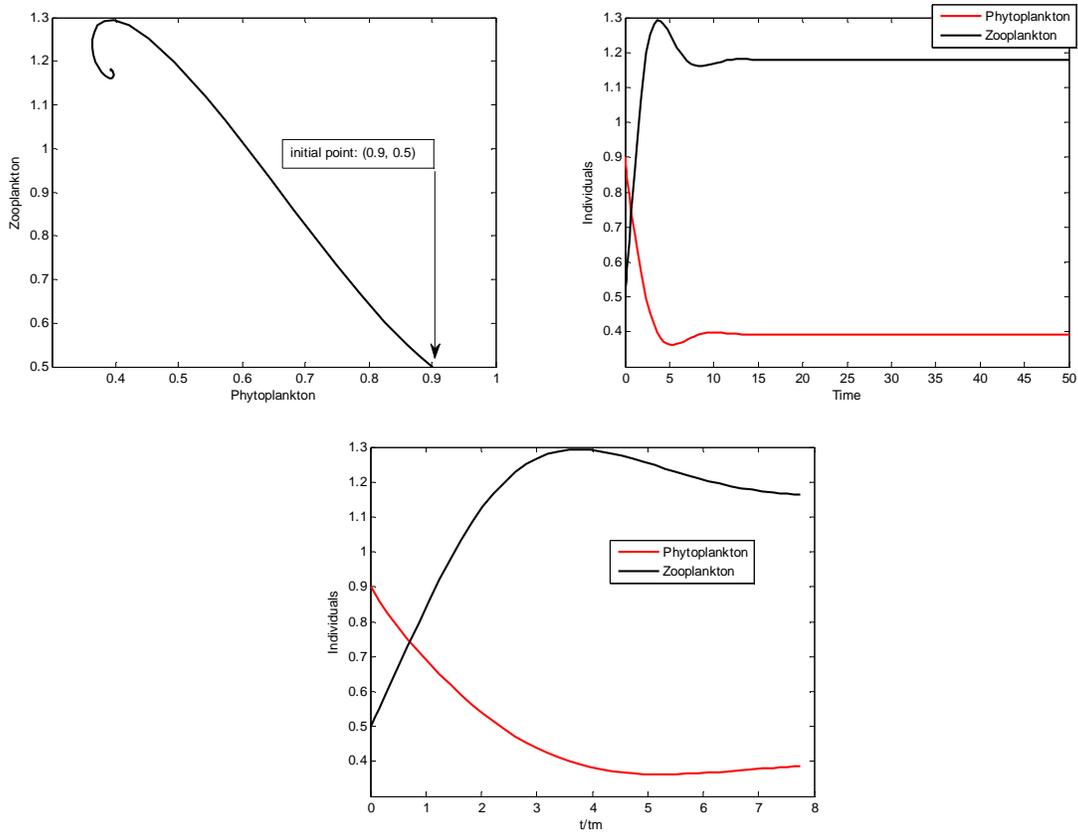

Fig. 1c. Phytoplankton-zooplankton evolution and orbit associated with $\mu = 1$.

For the parameters set, $\alpha = 1, \lambda = 0.057, \beta = 1.2, \gamma = 0.5, t_m = 8.88, x_0 = 0.5, y_0 = 0.5$, we have plotted the phase portrait and its time-series which are given below (2a, 2b, 2c):

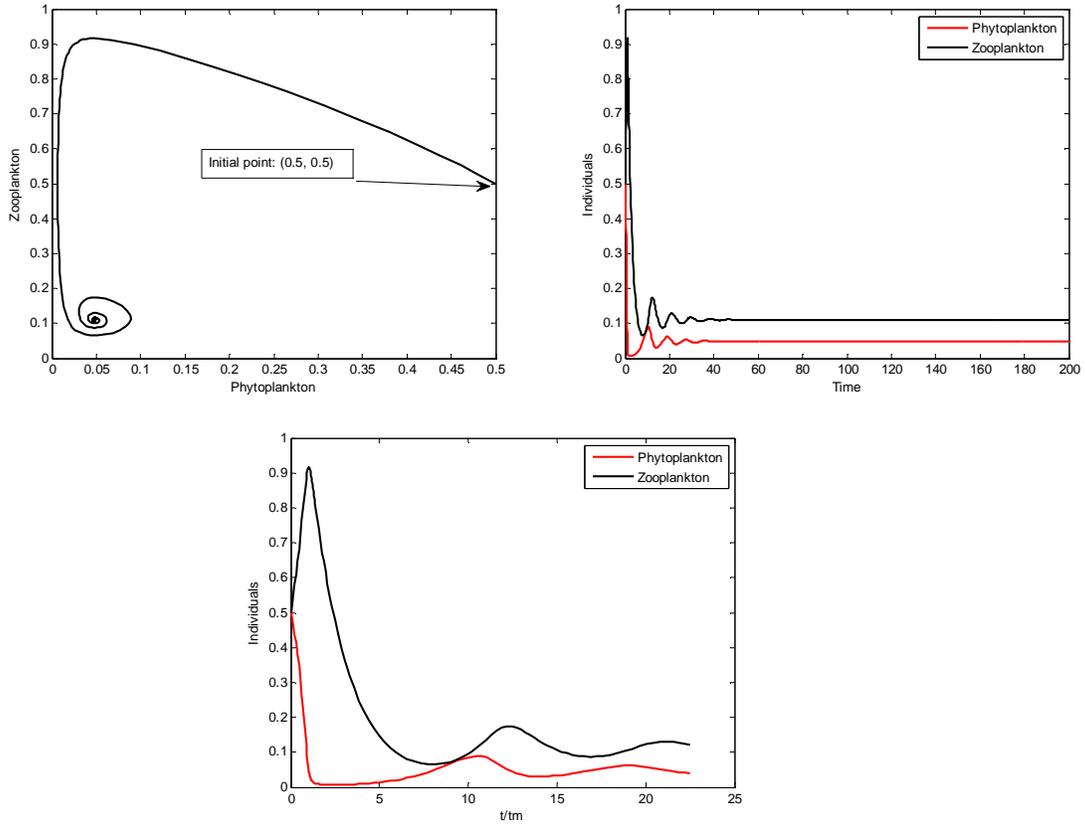

Fig. 2a. Phytoplankton-zooplankton evolution and orbit associated with $\mu = 0$.

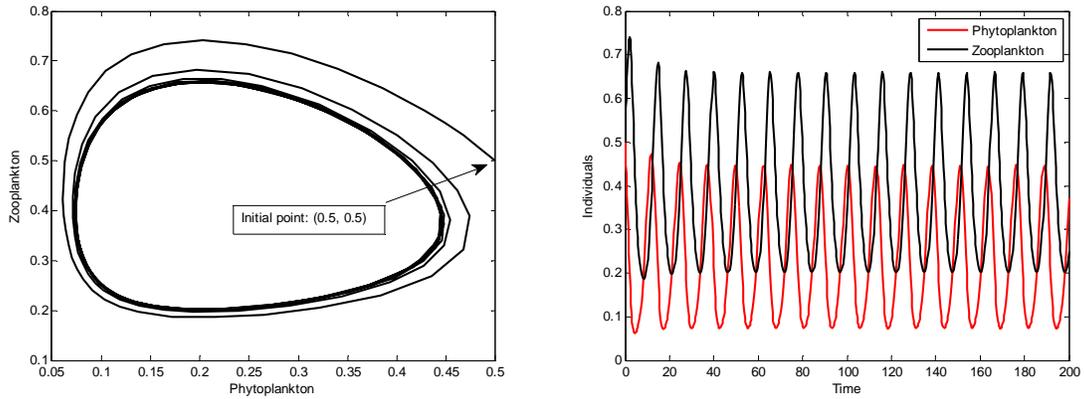

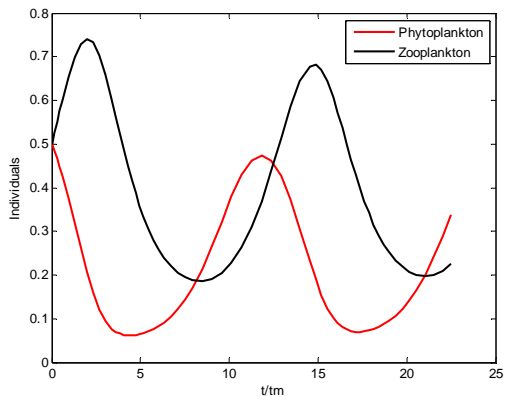

Fig. 2b. Phytoplankton-zooplankton evolution and orbit associated with $\mu = 0.5$.

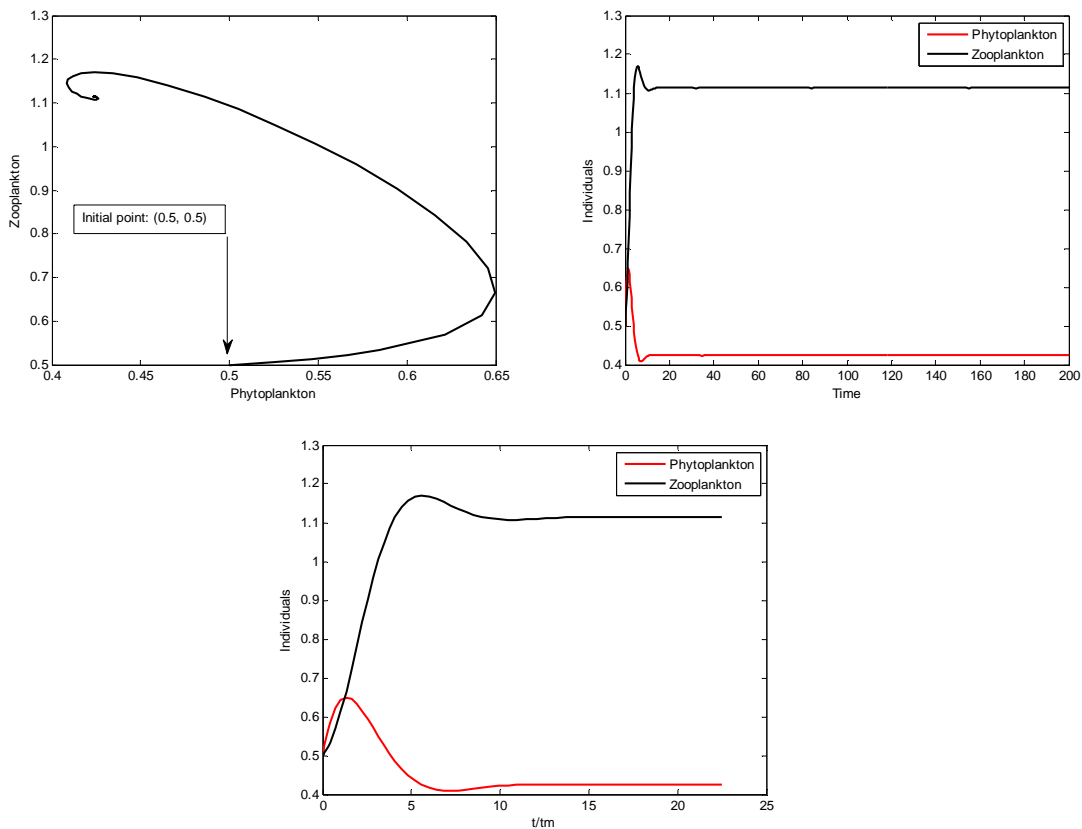

Fig. 2c. Phytoplankton-zooplankton evolution and orbit associated with $\mu = 1$.

## 5. Conclusions and Discussions

The similarity of the behavior of the two populations, with respect to the results presented by Freund *et al*. [18], Martin [19], and Malchow *et al*. [20], suggests that the results may be applicable to a broad class of zoological populations, bi-stable systems which exhibit excitable dynamics. Some authors [21-23] have studied the affects of variations in the initial composition of the system on the time dependent behavior. In this context, we have considered the affects of different initial conditions (1:1 and 1.8:1 proportions). Results show that the asymptotic behaviors of the analyzed models agree with those published by Richards *et al*. [23], and Biktashev *et al*. [24].

**References**

[1]  Truscott, J., Brindely, J., 1994. Ocean plankton populations as excitable media. Bulletin of Mathematical Biology, 56, No. 1, 185-211.

[2]  Holling, C. S. 1959. The components of predation as revealed by a study of small mammal predation on the European pine sawfly. Can. Ent. 91, 293.

[3]  Uye, S. 1986. Impact of copepod grazing on the red tide flagellate chattonella antique. Mar. Bio.92, 35.

[4]  Wake, G. et al. 1991. Picoplankton and marine food chain dynamics in a variable mixed layer: a reaction-diffusion model. Ecol. Modelling 57, 193.

[5]  Oyodum, O.D., Awojoyogbe, O.B., Dada, M., Magnuson, J., 2009. On the earliest definition of the Boubaker polynomials. Eur. Phys. J. Appl. Phys. 46, 21201–21203.

[6]  Awojoyogbe, O.B., Boubaker, K., 2008. A solution to Bloch NMR flow equations for the analysis of homodynamic functions of blood flow system using m- Boubaker polynomials. Curr. Appl. Phys. 9, 278–283.

[7]  Ghanouchi, J., Labiadh, H., Boubaker, K., 2008. An Attempt to solve the heat transfer equation in a model of pyrolysis spray using 4q-order Boubaker polynomials. Int. J. Heat Technol. 26, 49–53.

[8]  Slama, S., Bouhafs, M., Ben Mahmoud, K.B., Boubaker, A., 2008. Polynomials solution to heat equation for monitoring A3 point evolution during resistance spot welding. Int. J. Heat Technol. 26, 141–146.

[9]  Slama, S., Bessrour, J., Bouhafs, M., BenMahmoud, K. B., 2009a. Numerical distribution of temperature as a guide to investigation of melting point maximal front spatial evolution during resistance spot welding using Boubaker polynomials.Numer.HeatTransferPartA55,401–408.

[10]  Slama, S., Boubaker, K., Bessrour, J., Bouhafs, M., 2009b. Study of temperature 3D profile during weld heating phase using Boubaker polynomials expansion. Thermochim. Acta 482, 8–11.

[11]  Fridjine, S., Amlouk, M., 2009. A new parameter: an ABACUS for optimizing functional materials using the Boubaker polynomials expansion scheme. Mod. Phys. Lett. B 23, 2179–2182.

[12]  Tabatabaei, S., Zhao, T., Awojoyogbe, O., Moses, F., 2009. Cut-off cooling velocity profiling inside a keyhole model using the Boubaker polynomials expansion scheme. Heat Mass Transfer 45, 1247–1251.

[13]  Zhao, T.G., Wang, Y.X., Ben Mahmoud, K.B., 2008. Limit and uniqueness of the Boubaker–Zhao polynomials imaginary root sequence. Int. J. Math. Comput. 1, 13–16.

[14]  Belhadj, A., Onyango, O., Rozibaeva, N., 2009. Boubaker polynomials expansion scheme-related heat transfer investigation inside keyhole model. J. Thermo- phys. Heat Transfer 23, 639–640.

[15]  Ghrib, T., Boubaker, K., Bouhafs, M., 2008. Investigation of thermal diffusivity–microhardness correlation extended to surface-nitrured steel using Boubaker polynomials expansion. Mod. Phys. Lett. B 22, 2893–2907.

[16]  Guezmir, N., Ben Nasrallah, T., Boubaker, K., Amlouk, M., Belgacem, S., 2009. Optical modeling of compound CuInS2 using relative dielectric function approach and Boubaker polynomials expansion scheme BPES. J. Alloys Compd. 481, 543–548.

[17]  Dubey, B., Zhao, T.G., Jonsson, M., Rahmanov, H., 2010. A solution to the accelerated-predator-satiety Lotka–Volterra predator–prey problem using Boubaker polynomial expansion scheme. J. Theor. Biology 264, 154-160.

[18]  Freund, J. A., Mieruch, S., Scholze, B., Wiltshire, K., Feudel, U., 2006. Bloom dynamics in a seasonally forced phytoplankton–zooplankton model: Trigger mechanisms and timing effects, Ecological Complexity, 3, 2, 129-139.

[19]  Martin, A.P., 2003, Phytoplankton patchiness: the role of lateral stirring and mixing Progress In Oceanography, 57 (2) 125-174.

[20]  Malchow, H., Hilker, F.M., Sarkar, R.R., Brauer, K., 2005. Spatiotemporal patterns in an excitable plankton system with lysogenic viral infection, Mathematical and Computer Modeling, 42, 9-10.

[21]  Sieber, M., Malchow, H., Schimansky-Geier, L., 2007. Constructive effects of environmental noise in an excitable prey–predator plankton system with infected prey, Ecological Complexity, 4 (4) 223-233.

[22]  Petrovskii, S., Li, B. L., Malchow, H., 2003. Quantification of the Spatial Aspect of Chaotic Dynamics in Biological and Chemical Systems, Bulletin of Mathematical Biology, 65, (3) 425-446.

[23]  Richards, K. J., Brentnall, S. J., 2006. The impact of diffusion and stirring on the dynamics of interacting populations, Journal of Theoretical Biology, 238 (2) 340-347.

[24] V. N. Biktashev, J. Brindley, 2004. Phytoplankton blooms and fish recruitment rate: effects of spatial distribution, Bulletin of Mathematical Biology, 66 (2) 233-259.

[25] Cressman, R., 2010. CSS, NIS and dynamic stability for two-species behavioral models with continuous trait spaces. Journal of Theoretical Biology 262, 80.

[26] Hauzy, C., Gauduchon, M., Hulot, F.D., Loreau, M., 2010. Density-dependent dispersal and relative dispersal affect the stability of predator–prey metacommunities. Journal of Theoretical Biology 266, 458.

[27] Townsend, S.E., Haydon, D.T., Matthews, L., 2010. On the generality of stability–complexity relationships in Lotka–Volterra ecosystems. Journal of Theoretical Biology 267, 243.